\documentclass[11pt]{amsart}

\usepackage{amsmath,amssymb,amsthm}
\usepackage{graphicx,pxfonts}

\DeclareMathAlphabet{\mathbbold}{U}{bbold}{m}{n}

\def\k{\mathbbold{k}}

\DeclareSymbolFont{rsfscript}{OMS}{rsfs}{m}{n}
\DeclareSymbolFontAlphabet{\mathrsfs}{rsfscript}

\DeclareFontFamily{OMS}{rsfs}{\skewchar\font'177}
\DeclareFontShape{OMS}{rsfs}{m}{n}{%
      <5> rsfs5
      <6> <7> rsfs7
      <8> <9> <10> rsfs10
      <10.95> <12> <14.4> <17.28> <20.74> <24.88> rsfs10
      }{}

\def\calA{\mathrsfs{A}}
\def\calB{\mathrsfs{B}}

\def\calF{\mathrsfs{F}}
\def\calG{\mathrsfs{G}}

\def\calK{\mathrsfs{K}}

\def\calM{\mathrsfs{M}}
\def\calN{\mathrsfs{N}}
\def\calO{\mathrsfs{O}}
\def\calP{\mathrsfs{P}}
\def\calQ{\mathrsfs{Q}}
\def\calR{\mathrsfs{R}}
\def\calS{\mathrsfs{S}}

\def\calV{\mathrsfs{V}}
\def\calW{\mathrsfs{W}}

\DeclareMathOperator{\As}{As}
\DeclareMathOperator{\Mag}{Mag}
\DeclareMathOperator{\Lie}{Lie}
\DeclareMathOperator{\Com}{Com}

\DeclareMathOperator{\Perm}{Perm}
\DeclareMathOperator{\Dend}{Dend}
\DeclareMathOperator{\Dias}{Dias}

\DeclareMathOperator{\PreLie}{PreLie}

\DeclareMathOperator{\Fin}{Fin}
\DeclareMathOperator{\Ord}{Ord}
\DeclareMathOperator{\Vect}{Vect}

\makeatletter

\theoremstyle{plain}

\newtheorem {theorem}{Theorem}
\newtheorem {lemma}{Lemma}
\newtheorem {corollary}{Corollary}

\newtheorem {proposition}{Proposition}

\theoremstyle{definition}

\newtheorem {definition}{Definition}
\newtheorem {remark}{Remark}
\newtheorem {example}{Example}

\begin{document}
\title{Freeness theorems for operads via Gr\"obner bases}
\author{Vladimir Dotsenko}
\address{Dublin Institute for Advanced Studies, 10 Burlington Road, Dublin 4, Ireland and School of Mathematics, Trinity College, Dublin 2, Ireland}
\email{vdots@maths.tcd.ie}

\thanks{The author's research was supported by an IRCSET research fellowship.}

\begin{abstract}
We show how to use Gr\"obner bases for operads to prove various freeness theorems: freeness of certain operads as nonsymmetric operads,
freeness of an operad $\calQ$ as a $\calP$-module for an inclusion $\calP\hookrightarrow\calQ$, freeness of a suboperad. This gives new proofs of many 
known results of this type and helps to prove some new results. 
\end{abstract}

\maketitle

\section{Introduction}

\subsection{Description of results}
Recently, many freeness theorems about operads and free algebras over various operads have been proved. An incomplete list includes the following results:
\begin{itemize}
 \item free dendriform algebras are free as associative algebras (Loday and Ronco~\cite{LR});
 \item free pre-Lie algebras are free as Lie algebras (Chapoton~\cite{C} and Foissy~\cite{F-PreLie}); 
 \item free algebras with two compatible associative products are free as associative algebras (the author's result~\cite{D});
 \item the nonsymmetric operads $\Lie$ and $\PreLie$ are free (Salvatore and Tauraso~\cite{ST}, Bergeron and Livernet~\cite{BLi});
 \item the suboperad of the operad $\PreLie$ generated by the symmetrized pre-Lie product is free (Bergeron and Loday~\cite{BLo}).
\end{itemize}

In this article, we apply Gr\"obner bases for operads to derive several freeness theorems that imply all these results and several new ones. Our freeness theorems remind of the Magnus's Freiheitssatz from the group 
theory~\cite{Mag} and its analogues in other branches of algebra. Many of the results of this paper can be obtained by direct computations that use the Gr\"obner basis algorithm~\cite{DK}, however, we tried to replace most of 
computations by ideas coming from the Koszul duality theory~\cite{GK}.

Besides obvious applications of freeness theorems to computations of dimensions and bases for our operads, we hope that in some cases our results have further applications. In particular, results on freeness as a module are interesting from the homological algebra point of view: free modules are used to construct resolutions, so our freeness theorems may be used to ensure freeness of resolutions. We hope to address that elsewhere, together with various questions of shuffle homological algebra for non-free modules.

\subsection{Outline of the paper}

The paper is organized as follows. In Section~\ref{Groebner}, we briefly recall shuffle operads and Gr\"obner bases. In Section~\ref{NonSym}, we prove a general criterion for a symmetric operad to be free as a nonsymmetric operad, 
and show how this criterion applies to the cases of $\Lie$, $\PreLie$, and $\Lie^2$ (first two of these were establised earlier \cite{ST,BLi}, the last one is new). In Section~\ref{Suboperad}, we prove a criterion for a mapping of operads to be an embedding, and apply it to deduce the inclusion $\Mag\hookrightarrow\PreLie$ (proved by Bergeron and Loday by different methods in a forthcoming paper~\cite{BLo}). In Section~\ref{Module}, we prove a criterion of freeness as a module, and show how this criterion applies to the known cases $(\Lie,\PreLie)$~\cite{C,F-PreLie} and $(\As,\Dend)$~\cite{LR}, and, for a certain class of quadratic operads, to the case $(\calO,\calO^2)$ 
where $\calO^2$ denotes the operad of (weakly) compatible $\calO$-structures (generalizing our earlier result on compatible associative products~\cite{D}).

\subsection{Acknowledgements}
The author is grateful to Jean-Louis Loday and Bruno Vallette for an invitation to the conference ``Operads 2009'' where he came up with some of the ideas used in this article. He is also grateful to Paolo 
Salvatore and Muriel Livernet for some interesting discussions. Special thanks are due to Jean-Louis Loday for sharing an unpublished result on the symmetrized pre-Lie product.

\section{Shuffle operads and Gr\"obner bases}\label{Groebner}

All vector spaces throughout this work are defined over an arbitary field~$\k$ of zero characteristic. 

In this section, we give, mostly following \cite{DK}, a brief outline of definitions and the most important facts. For details on symmetric operads and Koszul duality, see \cite{GK} and \cite{MSS}. For more details on shuffle operads and Gr\"obner 
bases, see \cite{DK}.

\subsection{Shuffle compositions}

We denote by~$\Ord$ the category of nonempty finite ordered sets (with order-preserving bijections as morphisms), and by $\Fin$~--- the category of nonempty finite sets 
(with bijections as morphisms). Also, we denote by $\Vect$ the category of vector spaces (with linear operators as morphisms; unlike the first two cases, we do not require a map to be invertible).

\begin{definition}
\begin{enumerate}
 \item A \emph{(nonsymmetric) collection} is a functor from the category~$\Ord$ to the category~$\Vect$.
 \item A \emph{symmetric collection} (or an \emph{$\mathbb{S}$-module}) is a functor from the category~$\Fin$ to the category~$\Vect$.
\end{enumerate}
For either type of collections, we can consider the category whose objects are collections of this type (and morphisms are morphisms of the corresponding functors). 
The natural forgetful functor ${}^{f}\colon\Ord\to\Fin$, $I\mapsto I^f$ leads to a forgetful functor ${}^f$ from the category of symmetric collections to the category of 
nonsymmetric ones, $\calP^f(I):=\calP(I^f)$. 
\end{definition}

The following monoidal structures on our categories are important for the theory of operads.

\begin{definition}
Let $\calP$ and $\calQ$ be two nonsymmetric collections. Define their \emph{(nonsymmetric) composition} $\calP\circ\calQ$ by the formula
 $$
(\calP\circ\calQ)(I):=\bigoplus_{k}\calP(k)\otimes\left(\bigoplus_{f\colon I\twoheadrightarrow[k]}\calQ(f^{-1}(1))\otimes\ldots\otimes\calQ(f^{-1}(k))\right),
 $$
where the sum is taken over all non-decreasing surjections~$f$.

Let $\calP$ and $\calQ$ be two nonsymmetric collections. Define their \emph{shuffle composition} $\calP\circ_{sh}\calQ$ by the formula
 $$
(\calP\circ_{sh}\calQ)(I):=\bigoplus_{k}\calP(k)\otimes\left(\bigoplus_{f\colon I\twoheadrightarrow[k]}\calQ(f^{-1}(1))\otimes\ldots\otimes\calQ(f^{-1}(k))\right),
 $$
where the sum is taken over all shuffling surjections~$f$, that is surjections for which~$\min f^{-1}(i)<\min f^{-1}(j)$ whenever~$i<j$.

Let $\calP$ and $\calQ$ be two symmetric collections. Define their \emph{(symmetric) composition} $\calP\circ\calQ$ by the formula
 $$
(\calP\circ\calQ)(I):=\bigoplus_{k}\calP(k)\otimes_{\k S_k}\left(\bigoplus_{f\colon I\twoheadrightarrow[k]}\calQ(f^{-1}(1))\otimes\ldots\otimes\calQ(f^{-1}(k))\right),
 $$
where the sum is taken over all surjections~$f$.
\end{definition}

\begin{definition}
\begin{enumerate}
 \item A \emph{nonsymmetric operad} is a monoid in the category of nonsymmetric collections with the monoidal structure given by the nonsymmetric composition.
 \item A \emph{shuffle operad} is a monoid in the category of nonsymmetric collections with the monoidal structure given by the shuffle composition. 
 \item A \emph{symmetric operad} is a monoid in the category of symmetric collections with the monoidal structure given by the (symmetric) composition.
\end{enumerate}
\end{definition}

It turns out that the forgetful functor is a monoidal functor between the category of symmetric operads and the category of shuffle operads. Consequently, it turns out that to study 
various questions of linear algebra for operads, it is sufficient to forget the full symmetric structure because the shuffle structure already captures everything.
Further in this section, the word ``operad'' means a shuffle operad.

\subsection{Tree monomials, divisibility, and Gr\"obner bases}

We use the usual way to represent operadic elements by decorated rooted trees. A tree has (internal) vertices, directed edges, and inputs (leaves). For a tree whose leaves are labelled by an ordered set, its canonical planar 
representative is defined as follows. In general, an embedding of a (rooted) tree in the plane is determined by an ordering of inputs for each vertex (in terms of the planar structure, it is the ordering of inputs from the left to the right). Thus, to define a planar embedding, we should be able to compare two inputs of every vertex~$v$. To do so, we find the minimal leaves that one can reach from~$v$ via the corresponding inputs. The input for which the minimal leaf is smaller is considered to be less than the other one. Note that this choice of a representative is essentially the same one as we already made when we identified symmetric compositions with shuffle compositions.

Let us introduce an explicit realisation of the free operad generated by a collection $\calV$. The basis of this operad will be indexed by planar representative of trees with decorations of all vertices. First of all, the 
simplest possible tree is the degenerate tree (without internal vertices); it corresponds to the unit of our operad. The second simplest type of trees is given by corollas, that is trees with one vertex. We shall fix a 
basis~$B^\calV$ of $\calV$ and decorate the vertex of each corolla with a basis element; for a corolla with $n$ inputs, the corresponding element should belong to the basis of~$\calV(n)$. The basis for whole free operad consists 
of all planar representatives of trees built from these corollas (explicitly, one starts with this collection of corollas, defines compositions of trees in terms of grafting, and then considers all trees obtained from corollas by 
iterated shuffle compositions). We shall refer to elements of this basis as \emph{tree monomials}.

An ordering of tree monomials of $\calF_\calV$ is said to be \emph{admissible}, if it is compatible with the operadic structure, that is, replacing the operations in any shuffle compositions with larger operations of the same 
arities increases the result of the composition. Here we shall describe several admissible orderings which suit our purposes. All results of this section are valid for every admissible ordering of tree monomials.

Recall the following construction crucial for the ``path-lexicographic ordering''~\cite{DK}. Let $\alpha$ be a tree monomial with $n$ inputs. We associate to $\alpha$ a sequence $(a_1,a_2, \ldots,a_n)$ of $n$ words in the 
alphabet $B^\calV$ and a permutation $g\in S_n$ as follows. For each leaf~$i$ of the underlying tree $\tau$, there exists a unique path from the root to $i$. The word $a_i$ is the word composed, from left to right, of the labels 
of the vertices of this path, starting from the root vertex. The permutation $g$ lists the labels of leaves of the underlying tree in the order determined by the planar structure (from left to right).

To compare two tree monomials we always compare their arities first. If the arities are equal, there are several different options of how to proceed. Recall that an ordering of words in the alphabet $B^\calV$ is said to be 
admissible, if it is compatible with the semigroup structure on words (increasing factors increases the product). Fix some admissible ordering of words: that can be the lexicographic ordering, or its reverse, or the 
degree-lexicographic ordering (compare lengths first, and if they are equal, then compare lexicographically, or its reverse, or something else of the same sort). Now, sequences of words can be compared lexicographically: first 
compare the first words in both sequences, if they are equal, compare the second words etc. Permutations may be compared in the lexicographic order or its reverse. Also, the result depends on what we compare first, the 
permutations or the sequences of words. This gives rise to many candidates for an ordering. Generalizing~\cite{DK,H}, one can prove

\begin{proposition}
For each admissible ordering of words, all the orderings described above are admissible. 
\end{proposition}

For a tree monomial~$\alpha$ with the underlying tree $T$ and a subtree~$T'$ of~$T$, let us define a tree monomial~$\alpha'$ that corresponds to~$T'$. Its vertices are already decorated, so we just need to take care of 
the leaf labelling. For each leaf $l$ of $T'$, let us consider the smallest leaf of $T$ that can be reached from~$l$. We then number the leaves according to these ``smallest descendants'': the leaf with the smallest possible 
descendant gets the label~$1$, the second smallest~--- the label~$2$ etc.

\begin{definition}
For two tree monomials $\alpha$, $\beta$ in the free operad $\calF_\calV$, we say that \emph{$\alpha$ is divisible by $\beta$}, if there exists a subtree of the underlying tree of $\alpha$ for 
which the corresponding tree monomial $\alpha'$ is equal to~$\beta$.
\end{definition}

\begin{definition}
For an element~$\lambda$ of the free operad, the tree monomial~$\alpha$ is said to be its \emph{leading term}, if it is the largest of the terms which occur (with a nonzero coefficient) in the expansion of $\lambda$ as a linear combination of tree monomials. 
\end{definition}

Here and below we assume that $\calM$ is an operadic ideal of~$\calF_\calV$, and $\calG$ is a system of generators of~$\calM$. 

\begin{definition}
$\calG$ is called a \emph{Gr\"obner basis} of~$\calM$, if for every $f\in\calM$ the leading term of $f$ is divisible by the leading term of some element of~$\calG$.
The element $f\in \calF_\calV$ is said to have \emph{the residue $\overline{f}$ modulo $\calG$}, if $f-\overline{f}\in\calM$, and $\overline{f}$ is a linear combination of tree monomials that are not divisible by leading terms of 
elements of~$\calG$ (\emph{normal tree monomials}).  Notation: $f\equiv\overline{f}\pmod{\calG}$.

\end{definition}

For an explicit algorithm computing Gr\"obner bases, see \cite{DK}. Here we only use some simple criteria for Gr\"obner bases.

\begin{proposition}[\cite{DK}]
The set $\calG$ is a Gr\"obner basis for $\calM$ if and only if the (images of) normal tree monomials form a basis of the quotient~$\calF_\calV/\calM$.
\end{proposition}

\begin{definition}
Let $\calP$ be an operad, $\calP\simeq\calF_\calV/\calM$. A set of tree monomials $B^\calP\supset\calV$ in the free operad $\calF_\calV$ is said to be a \emph{$k$-triangular basis} of~$\calP$ if 
\begin{enumerate}
\item The image of $B^\calP$ under the canonical projection $$\calF_\calV\twoheadrightarrow \calF_\calV/\calM\simeq\calP$$ is a basis of $\calP$.
\item Every shuffle composition of tree monomials from~$B^\calP$ is either an element of~$B^\calP$ or is congruent modulo~$\calM$ to a linear combination of smaller elements of~$B^\calP$.
\item A tree monomial $\alpha$ belongs to $B^\calP$ if and only if for every its subtree with at most $k$ vertices the corresponding divisor belongs to $B^\calP$.	
\end{enumerate}
\end{definition}

\begin{proposition}[\cite{DK}]\label{k-triang}
Let $\calP\simeq\calF_\calV/\calM$ be an operad, $\calG\subset\calM$ be a system of generators. Then if $\calG$ is a Gr\"obner basis for $\calM$, then the set of tree monomials which are not divisible by the leading terms of 
elements of~$\calG$ is a $k$-triangular basis of $\calP$, where $k$ is the maximal number of vertices in those leading terms. Conversely, for a $k$-triangular basis~$\calB$, there exists a Gr\"obner basis whose elements 
are combinations of tree monomials with at most~$k$ vertices. In this case, $\calB$ coincides with the set of all normal tree monomials.
\end{proposition}

Denote by $LT_2(\calP)$ the set of the leading terms of the quadratic part of the Gr\"obner basis for the operad $\calP$ (that is, the elements of the Gr\"obner basis which are combinations of tree monomials with two internal vertices). Let us identify, on the level of vector spaces, the free operads whose quotients are, respectively, $\calP$ and its Koszul dual $\calP^!$ (identifying the dual bases of corollas, and the corresponding tree monomials), taking for the ordering of tree monomials for~$\calP^!$ the one opposite to the 
ordering for~$\calP$. In the case when we have elements with more than~$2$ inputs among the operad generators, working with the Koszul dual operad  brings us to the realm of dg-operads. In that case, the machinery of Gr\"obner bases exists, but when trying to compute something one should be careful: a tree monomial is defined up to a sign and this should be taken into account when collection terms in a result of a computation. 

The following can be easily derived from~\cite{H}.

\begin{proposition}
The set $LT_2(\calP^!)$ is the set of quadratic tree monomials complementary to~$LT_2(\calP)$. Also, $\calP$ has a quadratic Gr\"obner basis if and only if $\calP^!$ has a quadratic Gr\"obner basis; in this case both $\calP$ and $\calP^!$ are Koszul. 
\end{proposition}

\begin{corollary}
Under the conditions of the previous proposition, the set of tree monomials for which every quadratic divisor is a leading monomial of an element of the Gr\"obner basis of $\calP$ spans the Koszul dual~$\calP^!$; 
the number of such monomials of arity~$n$ gives an upper bound on~$\dim\calP^!(n)$. If this upper bound is sharp for all~$n$, $\calP^!$ has a quadratic Gr\"obner basis. 
\end{corollary}

\section{Freeness as a nonsymmetric operad}\label{NonSym}

The main idea of Salvatore and Tauraso in~\cite{ST} was to consider a special basis of the Lie operad, and prove that \emph{prime} elements of this basis (i.e. those which do not admit factorizations
as nonsymmetric compositions) generate this operad freely (as a nonsymmetric operad). We adapt their approach to prove a general freeness criterion.

\begin{definition}
Let $\alpha$ be a tree monomial. For each vertex of~$\alpha$, the (maximal) subtree rooted at this vertex is said to be \emph{connected} if the leaf labels of this subtree form an interval in the ordered set of leaves. A tree monomial is said to be \emph{prime} if its only vertex with a connected subtree is its root.
\end{definition}

\begin{example}
Among the tree monomials of arity~$3$
 $$
\includegraphics[scale=0.9]{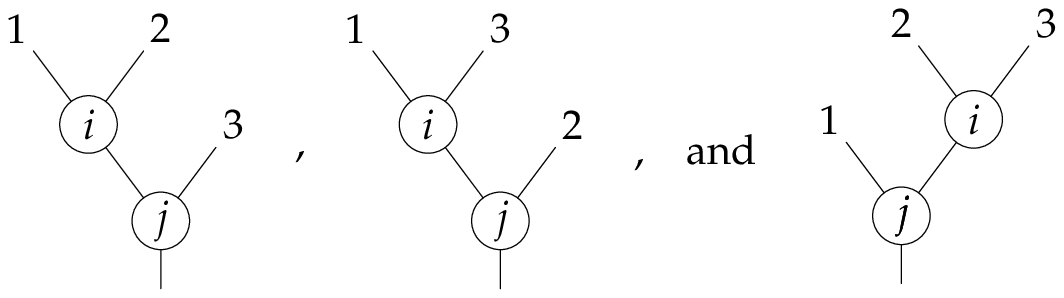}
 $$
monomials of the second type are prime, and others are not. 
\end{example}

\begin{theorem}\label{freeNS}
Let $\calO=\calF_\calV/(\calR)$ be a shuffle operad for which leading terms of a Gr\"obner basis consists of prime tree monomials. Then $\calO$ is free 
as a nonsymmetric operad.
\end{theorem}
\begin{proof}
Under the assumptions of the theorem, the nonsymmetric composition of two normal monomials is again normal, and a divisor of a normal monomial is normal. Thus, our statement follows immediately from 
\begin{lemma}
Let $\calO=\calF_\calV/(\calR)$ be a shuffle operad. Assume that there exists a family $\calB$ of monomials in $\calF_\calV$ satisfying the following conditions:
\begin{itemize}
 \item the images of monomials from $\calB$ under the natural projection form a basis of $\calO$; 
 \item $\calB$ is closed under nonsymmetric compositions, i.e. it forms a nonsymmetric suboperad of $\calF_\calV$;
 \item $\calB$ is closed under taking nonsymmetric divisors: if $\beta\in\calB$ is obtained from $\alpha$ by iterated nonsymmetric compositions, then $\alpha\in\calB$.
\end{itemize}
Then $\calO$ is free as a nonsymmetric operad. 
\end{lemma}
\begin{proof}
Consider the set $\calP$ of all prime monomials in $\calB$. We claim that the free nonsymmetric operad $\calF$ generated by $\calP$ is isomorphic to the nonsymmetric operad $\calB$ (which is isomorphic to $\calO$ as a 
nonsymmetric operad due to the first condition). Clearly, there exists a surjection $\phi\colon\calF\twoheadrightarrow\calB$. We shall describe an inverse mapping from $\calB$ to $\calF$. For a monomial $\gamma\in\calB$, consider 
the tree whose vertices correspond to vertices of $\gamma$ with connected subtrees; edges of this tree correspond to the partial order on vertices induced from the tree underlying $\gamma$. To such a vertex, we 
associate a prime monomial; namely, we take the unique prime nonsymmetric divisor rooted at the root of the original monomial (this divisor belongs to $\calB$ due to the third condition). Clearly, this mapping to $\calF$ is 
inverse to the surjection~$\phi$.
\end{proof}
\end{proof}

Before we formulate and prove the following result, we briefly recall several important definitions. The symmetric operad $\PreLie$ controlling pre-Lie algebras is generated by one operation $a,b\mapsto ab$ (without any 
symmetries) for which the associator $(a,b,c)=(ab)c-a(bc)$ is symmetric in the last two arguments. Its Koszul dual is denoted by $\Perm$, a basis of the space of $n$-ary operations of $\Perm$ is formed by commutative products of 
$n$ elements with one element emphasized. The symmetric operad $\Lie^2$ of compatible Lie brackets is generated by two skew-symmetric operations $a,b\mapsto [a,b]$ and $a,b\mapsto\{a,b\}$ which are compatible Lie brackets (every 
their linear combination is a Lie bracket). Its Koszul dual is denoted by ${}^2\Com$, it is generated by two symmetric binary operations which are totally compatible (every composition depends just on the types of the products 
used, not on how we compose them).

\begin{theorem}\label{NonSymFree}
The operads $\Lie$, $\PreLie$, and $\Lie^2$ are free as nonsymmetric operads.
\end{theorem}

\begin{proof}
Let us consider the ordering for which we first compare the permutations lexicographically, and in case when the permutations are equal~--- the sequences of words degree-lexicographically. 
We shall prove that in this case our operads have quadratic Gr\"obner bases with leading terms of the elements being prime monomials which, by Theorem~\ref{freeNS} is sufficient.

For the case of the operad $\Lie$, the leading term of its defining relation $[[a_1, a_2], a_3]-[[a_1, a_3], a_2]-[a_1, [a_2,a_3]]$ is $[[a_1, a_3], a_2]$. Note that this relation alone allows
the only normal tree monomial of arity~$n$ for the Koszul dual operad, namely 
 $$
[[[[a_1,a_n],a_{n-1}],\ldots,a_3],a_2],
 $$ 
so our upper bound on dimensions of components of the Koszul dual operad is sharp. This means that the Gr\"obner basis 
for the Koszul dual operad $\Com$ is quadratic, hence the Gr\"obner basis for $\Lie$ is quadratic. The leading term is indeed a prime monomial. 

For the case of the operad $\PreLie$, we use the operations $\alpha(a_1,a_2)=a_1a_2$ and $\beta(a_1,a_2)=a_2a_1$ which generate it as a shuffle operad. The leading terms of its defining relations
\begin{gather*}
\alpha(\alpha(a_1,a_2),a_3)-\alpha(a_1,\alpha(a_2,a_3))-\alpha(\alpha(a_1,a_3),a_2)+\alpha(a_1,\beta(a_2,a_3)),\\
\alpha(\beta(a_1,a_2),a_3)-\beta(\alpha(a_1,a_3),a_2)-\beta(a_1,\alpha(a_2,a_3))+\beta(\beta(a_1,a_3),a_2),\\
\alpha(\beta(a_1,a_3),a_2)-\beta(\alpha(a_1,a_2),a_3)-\beta(a_1,\beta(a_2,a_3))+\beta(\beta(a_1,a_2),a_3)
\end{gather*}
are (under the assumption $\beta>\alpha$) 
 $$
\alpha(\alpha(a_1,a_3),a_2), \beta(\beta(a_1,a_3),a_2), \text{ and } \alpha(\beta(a_1,a_3),a_2).
 $$
It follows that the normal monomials of arity~$n$ for the Koszul dual operad are among
 $$
\alpha(\ldots\alpha(\beta(\ldots\beta(a_1,a_n)\ldots),a_k),a_{k-1})\ldots),a_2),
 $$
with $i$ operations $\alpha$ and $j$ operations $\beta$, $i+j=n-1$. This upper bound for the Koszul dual operad $\Perm$ is sharp since $\dim\Perm(n)=n$. This means that there are no
more elements in the Gr\"obner basis and both $\Perm$ and $\PreLie$ have quadratic Gr\"obner bases. As we saw above, the leading terms of the relations for $\PreLie$ are prime monomials.

For the case of the operad $\Lie^2$, the leading terms of the defining relations
\begin{gather*}
[[a_1, a_2], a_3]-[[a_1, a_3], a_2]-[a_1, [a_2,a_3]],\\
\{\{a_1, a_2\}, a_3\}-\{\{a_1, a_3\}, a_2\}-\{a_1, \{a_2,a_3\}\},
\end{gather*}
and
\begin{multline*}
[\{a_1, a_2\}, a_3]-[\{a_1, a_3\}, a_2]-[a_1, \{a_2,a_3\}]+\\+
\{[a_1, a_2], a_3\}-\{[a_1, a_3], a_2\}-\{a_1, [a_2,a_3]\}
\end{multline*}
are (under the assumption $[\cdot,\cdot]>\{\cdot,\cdot\}$)
 $$
[[a_1, a_3], a_2], \{\{a_1, a_3\}, a_2\}, \text{ and } [\{a_1, a_3\}, a_2].
 $$
Similarly to the case of the operad $\PreLie$, this means that the dimension of the arity~$n$ component of the Koszul dual operad is at most $n$. Since the Koszul dual of $\Lie^2$ is ${}^2\Com$, and
$\dim(\ {}^2\Com(n))=n$, we deduce that both ${}^2\Com$ and $\Lie^2$ have quadratic Gr\"obner bases for our ordering. The leading terms of relations for $\Lie^2$ are prime monomials.
\end{proof}

\section{Embeddings of operads}\label{Suboperad}

\begin{proposition}\label{EmbOperads}
Let $\calP=\calF_\calV/(\calR)$ and $\calQ=\calF_{\calV\oplus\calW}/(\calR\oplus\calS)$ be shuffle operads. Assume that the Gr\"obner basis of $\calQ$ is a union of the Gr\"obner basis 
of $\calP$ and a set of relations whose leading terms are tree monomials not belonging to $\calF_\calV\subset\calF_{\calV\oplus\calW}$. Then the natural mapping $\calP\to\calQ$ is an embedding.
\end{proposition}

\begin{proof}
This statement is obvious: because of the condition imposed, the set of normal monomials for~$\calQ$ contains the set of normal monomials for~$\calP$.
\end{proof}

\begin{theorem}\label{MagPreLie}
The suboperad of the operad $\PreLie$ generated by the symmetrized pre-Lie product $a\cdot b=ab+ba$ is isomorphic to $\Mag$, the free operad on one generator. 
\end{theorem}

\begin{proof}
We consider the operad $\PreLie$ with its symmetric generator $\cdot$ and Lie bracket $[\cdot,\cdot]\colon a,b\mapsto [a,b]=ab-ba$. The following lemma is straightforward from the $\PreLie$ relations.
\begin{lemma}\label{PreLieSym}
The defining relations of $\PreLie$ for this choice of generators are
 $$
[[a_1, a_2], a_3]-[[a_1, a_3], a_2]-[a_1, [a_2,a_3]],
 $$
\begin{multline*}
 (a_1\cdot a_2)\cdot a_3-a_1\cdot(a_2\cdot a_3)-a_1\cdot[a_2,a_3]-[a_1,a_2]\cdot a_3-2[a_1,a_3]\cdot a_2+\\+[a_1,a_2\cdot a_3]+[a_1\cdot a_2,a_3]+[[a_1,a_3],a_2],
\end{multline*}
and
\begin{multline*}
 (a_1\cdot a_3)\cdot a_2-a_1\cdot(a_2\cdot a_3)+a_1\cdot[a_2,a_3]-[a_1,a_3]\cdot a_2-2[a_1,a_2]\cdot a_3+\\+[a_1,a_2\cdot a_3]+[a_1\cdot a_3,a_2]+[[a_1,a_2],a_3].
\end{multline*}
\end{lemma}
Let us choose an ordering of tree monomials as follows. To compare two tree monomials, we first compare their arities, then the number of the corollas corresponding to the product $\cdot$,
and then compare the tree monomials using the original path-lexicographic ordering from~\cite{DK}. It is easy to see that this ordering is admissible. 
The leading monomials of the relations of $\PreLie$ are (assuming that the bracket corolla is greater than the product corolla)
 $$
[[a_1, a_2], a_3], [a_1\cdot a_2,a_3], \text{ and } [a_1\cdot a_3,a_2].
 $$
It follows that the normal monomials of arity~$n$ for the Koszul dual operad are among
 $$
[[[[a_1\cdot a_k,a_n],a_{n-1}],\ldots,\hat{a}_k,\ldots,a_3],a_2], k\ge 2, \text{ and } [[[[a_1,a_n],a_{n-1}],\ldots,a_3],a_2],
 $$
which gives the upper bound $n$ on the dimensions of the arity~$n$ component of the Koszul dual operad, and, as we have seen before, this implies that our relations form a Gr\"obner basis of $\PreLie$ for this ordering. Clearly,
there are no leading terms of our relations that are made entirely of the product corollas, so by Proposition~\ref{EmbOperads} the product corolla defines an embedding $\Mag\hookrightarrow\PreLie$.
\end{proof}

\section{Freeness as a module}\label{Module}

\begin{theorem}\label{FreeModule}
Let $\calP\hookrightarrow\calQ$ be an embedding of shuffle operads, $\calP=\calF_\calV/(\calR)$, $\calQ=\calF_{\calV\oplus\calW}/(\calR\oplus\calS)$. Fix a basis of tree monomials in $\calF_{\calV\oplus\calW}$ for which the basis 
corollas is the union of the bases for~$\calV$ and~$\calW$.

\begin{enumerate}
\item Assume that the Gr\"obner basis of $\calQ$ is a union of the Gr\"obner basis of $\calP$ and a set of relations whose leading terms are tree monomials with the root from the basis of $\calW$. Then 
$\calQ$ is free as a left $\calP$-module, $\calQ\simeq\calP\circ\calK$ for some collection~$\calK$.
\item Assume that the Gr\"obner basis of $\calQ$ is a union of the Gr\"obner basis of $\calP$ and a set of relations whose leading terms are tree monomials where the parent vertex of each leaf belongs to basis of $\calW$. Then $\calQ$ is free as a right $\calP$-module, $\calQ\simeq\calK\circ\calP$ for some collection~$\calK$.
\end{enumerate}
\end{theorem}

\begin{proof}
These two statements are absolutely analogous to each other; let us prove the first one. Define $\calK$ to be the set of all normal monomials whose roots belong to the complementary set of generators~$\calW$. There exists a 
natural mapping $\calP\circ\calK\to\calQ$. This mapping is surjective for obvious reasons, moreover, it is injective because a composition of a normal tree monomial from $\calP$ and a normal tree monomial from $\calK$ cannot be reduced due to our 
condition on the Gr\"obner basis.
\end{proof}

\begin{corollary}
The operad $\PreLie$ is free as a $\Lie$-module. Consequently, free pre-Lie algebras are free as Lie algebras. 
\end{corollary}

\begin{proof}
Let us choose an ordering of tree monomials as follows. To compare two tree monomials, we first compare their arities, then corresponding permutations lexicographically, 
then the corresponding sequences of words degree-lexicographically.
The leading monomials of the relations of $\PreLie$ from Lemma~\ref{PreLieSym} are (assuming that the bracket corolla is greater than the product corolla)
 $$
[[a_1, a_3], a_2], [a_1,a_3]\cdot a_2, \text{ and } (a_1\cdot a_3)\cdot a_2.
 $$
Similarly to the proof of Theorem~\ref{NonSymFree}, this means that the dimension of the arity~$n$ component of the Koszul dual operad is at most $n$. Since the dimension of that component is equal to~$n$, we conclude that both 
$\PreLie$ and its dual have quadratic Gr\"obner bases for this ordering as well. For the last two relations which are complementary to the $\Lie$ relation, the roots of the leading terms belong to the complementary set of 
generators, which proves our statement.
\end{proof}

\begin{remark}
Theorem \ref{FreeModule}, applied to the Gr\"obner basis for the operad $\PreLie$ that we computed in Theorem~\ref{MagPreLie}, implies that $\PreLie$ is free as a left $\Mag$-module as well, $\PreLie\simeq\Mag\circ\calK$ for some 
collection~$\calK$. A direct computation shows that $\dim\calK(n)=n^{n-2}$, the number of labelled trees on~$[n]$. It would be interesting to describe $\calK$ combinatorially in the spirit of the description of $\PreLie$ via 
labelled rooted trees from~\cite{CL}.
\end{remark}

To formulate the main new result of this section, we recall the following general definitions that were given in~\cite{S}, as a generalisation of the notion of compatible Lie brackets (and the algebraic structure controlled by 
its Koszul dual operad) to general algebraic structures.

Let $\calO=\calF_\calV/(\calR)$ be a binary quadratic operad generated by binary operations $\alpha_1,\ldots, \alpha_s$ with $t$ relations
 $$
\left\{
\sum_{1\leq i,j \leq s} \gamma^{k,1}_{i,j}\alpha_j(\alpha_i(a_1,a_2),a_3)+\gamma^{k,2}_{i,j}\alpha_j(a_1,\alpha_i(a_2,a_3))
+\gamma^{k,3}_{i,j}\alpha_j(\alpha_i(a_1,a_3),a_2) \right\}_{1\leq k \leq t}.
 $$
Consider two operads $\calO_\circ=\calF_{\calV_\circ}/(\calR_\circ)$ and $\calO_\bullet=\calF_{\calV_\bullet}/(\calR_\bullet)$ both isomorphic to~$\calO$. We choose $\k$-bases $\alpha^1_1,\ldots,\alpha^1_s$ of $\calV_\circ$ and 
$\alpha^2_1,\ldots,\alpha^2_s$ of $\calV_\bullet$. The relations $\calR_\circ$ and $\calR_\bullet$ can then be given by the same $\gamma^{k,l}_{i,j}$. The suboperads $\calO_\circ$ and $\calO_\bullet$ of 
$\calF_{\calV_\circ\oplus\calV_\bullet}/(\calR_\circ\cup\calR_\bullet)$ control a pair of $\calO$-structures which are not related in any way. Now we are going to define two important ways of imposing additional compatibility 
relations.

\begin{definition}
Relations 
\begin{multline*}
\left\{\sum_{1\leq i,j \leq s} \gamma^{k,1}_{i,j}\alpha^1_j(\alpha^2_i(a_1,a_2),a_3)+\gamma^{k,2}_{i,j}\alpha^1_j(a_1,\alpha^2_i(a_2,a_3))
+\gamma^{k,3}_{i,j}\alpha^1_j(\alpha^2_i(a_1,a_3),a_2)+\right.\\ 
\left.\vphantom{\sum_{1\leq i,j \leq s}}+\gamma^{k,1}_{i,j}\alpha^2_j(\alpha^1_i(a_1,a_2),a_3)+\gamma^{k,2}_{i,j}\alpha^2_j(a_1,\alpha^1_i(a_2,a_3))+
\gamma^{k,3}_{i,j}\alpha^2_j(\alpha^1_i(a_1,a_3),a_2) \right\}_{1\leq k \leq t}
\end{multline*}
are called the linear compatibility (weak compatibility) relations. 
\end{definition}

\begin{definition}
Relations 
\begin{gather*}
 \alpha^1_j(\alpha^2_i(a_1,a_2),a_3)=\alpha^2_j(\alpha^1_i(a_1,a_2),a_3)\quad (1\le i,j\le s),\\
 \alpha^1_j(\alpha^2_i(a_1,a_3),a_2)=\alpha^2_j(\alpha^1_i(a_1,a_3),a_2)\quad (1\le i,j\le s),\\
 \alpha^1_j(a_1,\alpha^2_i(a_2,a_3))=\alpha^2_j(a_1,\alpha^1_i(a_2,a_3))\quad (1\le i,j\le s),\\
\left\{\sum_{1\leq i,j \leq s} \gamma^{k,1}_{i,j}\alpha^1_j(\alpha^2_i(a_1,a_2),a_3)+\gamma^{k,2}_{i,j}\alpha^1_j(a_1,\alpha^2_i(a_2,a_3))+
\gamma^{k,3}_{i,j}\alpha^1_j(\alpha^2_i(a_1,a_3),a_2) \right\}_{1\leq k \leq t},\\
\left\{\sum_{1\leq i,j \leq s} \gamma^{k,1}_{i,j}\alpha^2_j(\alpha^1_i(a_1,a_2),a_3)+\gamma^{k,2}_{i,j}\alpha^2_j(a_1,\alpha^1_i(a_2,a_3))+
\gamma^{k,3}_{i,j}\alpha^2_j(\alpha^1_i(a_1,a_3),a_2) \right\}_{1\leq k \leq t}
\end{gather*}
are called the total compatibility (strong compatibility) relations (the last one clearly follows from the first four and so is not necessary to include).
\end{definition}

\begin{definition}
\begin{enumerate}
 \item The operad $\calO^2$ is the quotient of $\calF_{\calV_\circ\oplus\calV_\bullet}$ modulo the ideal generated by $\calR_\circ$, $\calR_\bullet$, and the linear compatibility relations. It is called the operad of linearly 
compatible $\calO$-stuctures.
 \item The operad ${}^2\calO$ is the quotient of $\calF_{\calV_\circ\oplus\calV_\bullet}$ modulo the ideal generated by $\calR_\circ$, $\calR_\bullet$ and the total compatibility relations. It is called the operad of totally 
compatible $\calO$-structures.
\end{enumerate}
\end{definition}

This notation agrees with the notation introduced earlier for the particular cases $\Lie^2$ and $\,{}^2\Com$. 

\begin{proposition}[\cite{S}]
We have $(\calO^2)^!={}^2(\calO^!)$ and ${}^2\calO\simeq\,{}^2\Com\otimes\calO$.
\end{proposition}

\begin{theorem}\label{O2isPBW}
Let $\calO$ be an operad which has a quadratic Gr\"obner basis (for some admissible ordering). Furthermore, assume 
that all the monomials $\alpha_j(a_1,\alpha_i(a_2,a_3))$ are normal monomials relative to that Gr\"obner basis.
Then there exists an admissible ordering for which the operad $\calO^2$ has a quadratic Gr\"obner basis as well, and $\calO^2$ is a free $\calO$-module. 
In particular, free $\calO^2$-algebras are free as $\calO$-algebras.
\end{theorem}

\begin{proof}
Let us choose an ordering of tree monomials in $\calF_{\calV_\circ\oplus\calV_\bullet}$ as follows. To compare two tree monomials, we first compare their arities, then the number of corollas from the set of 
generators~$\calV_\circ$, then compare the images of these monomials under the natural homomorphism $\calF_{\calV_\circ\oplus\calV_\bullet}\twoheadrightarrow\calF_{\calV}$ (defined on generators by $\alpha^i_j\mapsto \alpha_j$), 
and if all the above procedures give the same result, compare these monomials using the path-lexicographic order, assuming that $\alpha^1_u>\alpha^2_v$ for all~$u,v$. It is clear that this ordering is admissible.

Note that from our normal monomial condition it follows that for the operad $\calO^!$ all the operations $\alpha_j(a_1,\alpha_i(a_2,a_3))$ are among the leading terms of the Gr\"obner basis. 
It follows that all underlying trees of normal tree monomials for this operad are ``combs'': all their vertices form a chain between the root and the leaf labelled~$1$. 

Let us consider the operad ${}^2(\calO^!)$. Since this operad is isomorphic to ${}^2\Com\otimes\calO^!$, we can easily construct a monomial basis for ${}^2(\calO^!)$ from a monomial basis of~$\calO^!$. Let us take the set of all 
normal tree monomials for $\calO^!$ and replicate every normal monomial of arity $n$ exactly $n$ times; for the $k^\text{th}$ replica ($1\le k\le n$), we look at the path from the root to the leaf labelled~$1$, and replace the last 
$k-1$ vertices of this path $\alpha_v$ by the corresponding generator $\alpha^1_v$ and the remaining $n-k$ vertices by the corresponding generator~$\alpha^2_v$. As a result, the corresponding tree monomial is the least (with 
respect to our ordering) among all the monomials that can be obtained this way (note that in ${}^2(\calO^!)$ all these monomials are equal to the least one because of the strong compatibility relations).

From Proposition~\ref{k-triang}, we see that normal tree monomials for $\calO^!$ form a $2$-triangular basis. This implies that our monomial basis for ${}^2(\calO^!)$ is $2$-triangular as well (the fact that all normal monomials 
are combs guarantees ``$2$'' in ``$2$-triangular''), so the total compatibility relations form a Gr\"obner basis. Consequently, the operad $\calO^2$ has a Gr\"obner basis as well, as the Koszul dual to~${}^2(\calO^!)$.

Finally, let us notice that for the leading terms of relations of $\calO^2$ which are not the relations of the first copy of $\calO_\circ$, the roots belong to~$\calV_\bullet$ (consistent orderings of monomials for an operad and 
its dual are opposite to each other), which proves the desired result.
\end{proof}

\begin{corollary}\label{O2isKoszul}
Under the assumption of Theorem~\ref{O2isPBW}, the operads~$\calO^2$ and ${}^2(\calO^!)$ are Koszul. 
\end{corollary}

Our result, despite of the additional assumption on the operad~$\calO$, are applicable to some new (not covered by~\cite{S}) cases, for example, the operad~$\Com$ (its Koszul dual~$\Lie$ is not set-theoretic, and so Strohmayer's 
results do not work, but our assumption applies). It turns out that some restrictions on the operad~$\calO$ should be imposed anyway (both for the Koszulness and for freeness of~$\calO^2$ as an $\calO$-module). Indeed, let us 
take the nilpotent free operad~$\calN$ with one symmetric generator. For this operad, it is easy to see that $\dim\calN^2(2)=2$, $\dim\calN^2(3)=3$, $\dim\calN^2(k)=0$ for~$k>3$. Thus the Hilbert series~\cite{GK} of this operad 
is $f(t)=t+t^2+\frac{t^3}{2}$, which immediately implies that this operad is not Koszul (since the inverse of~$f(-t)$ has negative coefficients). Also, we immediately see that~$\calN^2$ is not free as a left $\calN$-module (if it 
was, it would have an additional generator in arity~$2$ which, in turn, would lead to non-zero operations of arity~$4$).

It turns out that the result of Corollary~\ref{O2isKoszul} can be substantially generalized. For a binary quadratic operad $\calO=\calF_\calV/(\calR)$ let $\pi_\calO\colon\calF_\calV\to\calO$ be the natural projection. Let 
$\lambda$ be a tree monomial with $n-1$ vertices from the free operad on one generator. Fix an ordering of vertices of~$\lambda$, and let $\mathcal{L}^\calV_\lambda$ denote the induced decoration morphism from $\calV^{\otimes 
(n-1)}$ to $\calF_\calV$; it decorates vertices of $\lambda$ with elements of $V$. An operad~$\calO$ is said to be \emph{small} if for every $n\ge3$ and every tree monomial $\gamma$ with $n-1$ vertices the composite map 
$\pi_\calO\circ\mathcal{L}_\gamma^\calV\colon\calV^{\otimes (n-1)}\to\calO(n)$ is surjective (the class of small operads was studied in~\cite{Vallette} in relation to Manin products for operads; by \cite[Prop.~15]{Vallette}, 
smallness of a binary quadratic operad $\calA$ guarantees that for every binary quadratic operad~$\calB$ the Hadamard product~$\calA\otimes\calB$ and the white product~$\calA\medcirc\calB$ are isomorphic). Modelling the proof of 
Theorem~\ref{O2isPBW}, one can easily derive

\begin{proposition}
Suppose that the operad~$\calA$ is small. Under the assumption of Theorem~\ref{O2isPBW}, the operads~$\calO\medbullet\calA^!$ and $\calO^!\medcirc\calA$ have quadratic Gr\"obner bases and hence are Koszul.
\end{proposition}

\begin{remark}
In fact, in most of the cases of Koszul operads, we know an admissible ordering of tree monomials for which the Gr\"obner basis is quadratic. Moreover, we do not know any example of a Koszul operad generated by 
binary operations which has no quadratic Gr\"obner basis. Such examples should exist; in the case of quadratic algebras (that is, quadratic operads generated by unary operations), there are some known examples of that sort, 
see~\cite{Ber}.
\end{remark}

In the last example we wish to consider, the operads are regular, so everything is essentially reduced to the level of nonsymmetric operads. The symmetric operad $\Dend$ of dendriform algebras is generated by two binary 
operations $\prec$ and $\succ$ which satisfy the relations
\begin{gather*}
(a_1\prec a_2)\prec a_3=a_1\prec(a_2\prec a_3)+a_1\prec(a_2\succ a_3),\\
(a_1\succ a_2)\prec a_3=a_1\succ (a_2\prec a_3),\\
(a_1\prec a_2)\succ a_3+(a_1\succ a_2)\succ a_3=a_1\succ(a_2\succ a_3).
\end{gather*}
Its Koszul dual operad $\Dias$ of diassociative algebras is generated by two binary operations $\dashv$ and $\vdash$ which satisfy the relations
\begin{gather*}
(a_1\dashv a_2)\dashv a_3=a_1\dashv (a_2\dashv a_3),\quad (a_1\vdash a_2)\vdash a_3=a_1\vdash (a_2\vdash a_3),\\
(a_1\dashv a_2)\dashv a_3=a_1\dashv (a_2\vdash a_3),\quad (a_1\dashv a_2)\vdash a_3=a_1\vdash (a_2\vdash a_3),\\ 
(a_1\vdash a_2)\dashv a_3=a_1\vdash (a_2\dashv a_3).
\end{gather*}
It is well known that the operad $\Dias$ is a symmetrization of a nonsymmetric operad whose $n^\text{th}$ component has the dimension~$n$. We shall see that this allows for proofs similar to those for the pair of Koszul dual operads~$\PreLie$ and~$\Perm$. 

If we introduce the new generator $\star$, $a_1\star a_2=a_1\prec a_2+a_1\succ a_2$, it is easy to see that this operation is an associative product. In \cite{LR}, Loday and Ronco proved that the free dendriform algebra on one generator is free as an associative algebra with respect to this product. Since the dendriform operad is a symmetrization of a nonsymmetric operad, 
this essentially means that $\Dend$ is a free left module over $\As$ and hence the free dendriform algebra on any number of generators is free as an associative algebra. We shall prove 
this statement using Gr\"obner bases for nonsymmetric operads.

\begin{proposition}
The nonsymmetric operad $\Dend$ is a free left module over the nonsymmetric operad $\As$.
\end{proposition}

\begin{proof}
In terms of the generators $\star$ and $\succ$, the defining relations of $\Dend$ are
\begin{gather*}
(a_1\star a_2)\star a_3=a_1\star(a_2\star a_3),\\
(a_1\succ a_2)\succ a_3-(a_1\succ a_2)\star a_3- a_1\succ (a_2\succ a_3)+a_1\succ(a_2\star a_3)=0,\\
(a_1\star a_2)\succ a_3=a_1\succ(a_2\succ a_3).
\end{gather*}
If we assume that the operation $\succ$ is greater than the operation~$\star$, and use the path-lexicographic ordering of trees, the leading terms of relations are
$(a_1\star a_2)\star a_3$, $(a_1\succ a_2)\succ a_3$, and $(a_1\star a_2)\succ a_3$. We already know that this means an upper bound $n$ on the dimension of the $n^\text{th}$ component of 
the Koszul dual operad. Thus, the Gr\"obner basis is quadratic since this bound is sharp. The leading terms of the relations (except for the associativity relation) have their roots labelled by
$\succ$ which proves that $\Dend$ is a free left $\As$-module.
\end{proof}

\begin{remark}
A freeness theorem which is not covered by this text is the result of Foissy~\cite{F-Brace} stating that free brace algebras are free as pre-Lie algebras. 
It remains an open question to compute the Gr\"obner basis for the brace operad and derive the result of Foissy from Theorem~\ref{FreeModule}.
\end{remark}

\end{document}